\documentclass[11pt, a4paper]{amsart}
\usepackage{amsfonts,amsthm,amssymb,amsmath}

\numberwithin{equation}{section}
\usepackage{latexsym}
\usepackage[dvips]{graphicx}
\usepackage{geometry}
\title {Entropies of compact strictly convex projective manifolds\\}
\author{Micka\"el Crampon}
\email{crampon@math.u-strasbg.fr}

\newtheorem{defi}{Definition}[section]

\newtheorem{thm}[defi]{Theorem}
\newtheorem*{thmsansnum}{Theorem}
\newtheorem{prop}[defi]{Proposition}
\newtheorem{lemma}[defi]{Lemma}
\newtheorem{cor}[defi]{Corollary}

\def\C{\mathbb{C}}

\def\P{\mathbb{P}}

\def\R{\mathbb{R}}

\def\ph{\varphi}

\def\ho{H\Omega}
\def\o{\Omega}
\def\g{\gamma}
\def\l{\lambda}
\def\ch{\chi}
\def\M{\mathcal{M}}
\def\vect{\textrm{vect}}
\def\G{\mathcal{G}}

\begin{document}

\begin{abstract}
Let $M$ be a compact manifold of dimension $n$ with a strictly convex projective structure. We consider the geodesic flow of the Hilbert metric on it which is known to be Anosov. We prove that its topological entropy is less than $n-1$, with equality if and only if the structure is Riemannian, that is hyperbolic. As a corollary, we get that the volume entropy of a divisible strictly convex set is less than $n-1$, with equality if and only if it is an ellipsoid.
\end{abstract}

\maketitle

\section{Introduction}

In 1936, in what seems to be the first general introduction to the notion of locally homogeneous space \cite{ehresmann}, Charles Ehresmann noticed the following : it is not excluded for the universal covering of some compact locally projective surface to be a bounded convex domain whose boundary would not be analytic. But immediately he added that to his mind such a case should not occur.
Thirty years later, Kac and Vinberg \cite{kav} proved this implausible situation was indeed possible.\\
Such surfaces, and by extension, such manifolds are the main objects of this article. These are compact manifolds which can be written as a quotient $\o/\Gamma$, where $\o$ is a strictly convex proper open set of the projective space and $\Gamma$ a subgroup of the projective group acting cocompactly on $\o$. Such a manifold is said to be strictly convex projective and $\o$ is said to be divisible. From now, all the manifolds are compact.\\\\
Lots of manifolds admit strictly convex projective structures. The basic example is a hyperbolic manifold for which the Beltrami-Klein model of hyperbolic space provides such a structure. As observed in many occasions by various authors, any other strictly convex projective structure is much less regular, that is the boundary of the convex set is so. Ehresmann noticed first that it was nowhere analytic. Then Benz\'ecri \cite{benzecri} proved that if the boundary was $C^2$, then the convex set was an ellipsoid. Finally, from a different point of view, Edith Soci\'e-Methou \cite{sm} proved that the ellipsoid was the only convex set whose boundary is $C^2$ and with positive definite Hessian, on which a subgroup of the projective group could act properly.\\
Despite everything, these non-smooth structures are numerous in various senses :
\begin{itemize}
\item If a manifold admits a hyperbolic structure then it also admits non-smooth strictly convex projective structures ; more over, the deformation space $\G (M)$ of such structures is really bigger than the Teichm\"uller space $\mathcal{T}(M)$ of hyperbolic structures.\\
In dimension 2, Goldman \cite{goldman} proved that it was an open cell of dimension $8g-8$, where $g\geq 2$ denotes the genus of the surface, whereas $\mathcal{T}(M)$ is only of dimension $3g-3$.\\
In dimension higher than $3$, Mostow's rigidity theorem \cite{mostow} affirms that $\mathcal{T}(M)$ is reduced to a point. But it follows from the works of Benz\'ecri \cite{benzecri} and Koszul \cite{koszul} on affine and projective manifolds that $\G(M)$ is open in the space of projective structures $\R\P^n(M)$. In particular, Johnson and Millson \cite{johnsonmillson} constructed non trivial continuous deformations of a hyperbolic structure into strictly convex projective ones.
\item There are manifolds which admit strictly convex projective structures but no hyperbolic structure. Such example cannot exist in dimension 2 and 3 but Benoist \cite{benoistqi} constructed an example in dimension 4, and Kapovich \cite{kapo} proved that some Gromov-Thurston manifolds \cite{gromovthurston} actually provided other examples.\\
\end{itemize}
Any strictly convex set $\o$ carries a Hilbert metric $d_{\o}$ (see section \ref{hilbertgeom}). When $\o$ is an ellipsoid, $(\o,d_{\o})$ coincides with the hyperbolic space; in the other cases, the metric is not Riemannian anymore, but comes from a Finsler metric which has the same regularity as the boundary of the convex. Hilbert metric is invariant under any homography, and thus provides a metric on any compact projective manifold $M=\o/\Gamma$. With this metric, $M$ is projectively flat : in local projective charts, geodesics, as locally shortest paths, are straight lines.\\\\
These structures are for various reasons generalizations of hyperbolic ones. Despite the lack of regularity, we can define a notion of curvature and prove it is constant and strictly negative. Furthermore, Yves Benoist proved the following theorem :
\begin{thmsansnum}[\cite{benoistcv1}]
Let $\o$ be a divisible convex set, divided by $\Gamma$. The following statements are equivalent :
\begin{itemize}
\item the space $(\o,d_{\o})$ is Gromov-hyperbolic ;
\item $\o$ is strictly convex ;
\item the boundary $\partial\o$ of $\o$ is $C^1$ ;
\item $\Gamma$ is Gromov-hyperbolic.\\
\end{itemize}
\end{thmsansnum}
This paper can be seen as a continuation of \cite{benoistcv1}, where Benoist initiated the study of the geodesic flow of the Hilbert metric. In particular, Benoist proved similar properties to those of the hyperbolic geodesic flow, namely that the flow was Anosov and topologically mixing. But he already made the following observation, which distinguished the two dynamical systems : whereas hyperbolic geodesic flows admit the Liouville measure as natural invariant measure, the others do not admit any smooth invariant measure.\\
A major invariant in the theory of dynamical systems (see \cite{hist}) is the topological entropy, which roughly speaking measures how the system separates the points, how much it is chaotic. Let us recall briefly its definition. Given a system $\ph^t:X\longrightarrow X$ on a compact metric space $(X,d)$, we define the distances $d_t,\ t\geq 0$, on $X$ by $d_t(x,y)=\max_{0\leq s\leq t}\ d(\varphi^s(x),\varphi^s(y)),\ x,y\in X.$ The topological entropy of $\ph$ is then the well defined quantity
$$h_{top}(\varphi)=\lim_{\epsilon \to 0}\Big[ \limsup_{t \to \infty}\frac{1}{t}\log N(\varphi,t,\epsilon)\Big] \in [0,+\infty],$$
where $N(\varphi,t,\epsilon)$ denotes the minimal number of open sets of diameter less than $\epsilon$ for $d_t$ needed to cover $X$.\\
It is well known that the topological entropy of the hyperbolic geodesic flow is $n-1$ when the manifold is of dimension $n$. Our main theorem answers a question that raised during a Finsler meeting at the CIRM in 2005 and provides a new distinction between the Riemannian and the non-Riemannian cases :

\begin{thm}\label{majeur}
Let $\ph$ be the geodesic flow of the Hilbert metric on a strictly convex projective compact manifold $M$ of dimension $n$. Its topological entropy $h_{top}(\ph)$ satisfies the inequality
$$h_{top}(\ph) \leq (n-1),$$ 
with equality if and only if the Hilbert metric comes from a Riemannian metric.
\end{thm}

The proof of this result is mainly based on results in the Anosov systems theory, developed since the 60's, and on the geometrical approach to second order differential equations made by Patrick Foulon in \cite{foulon86}. (See also the appendix of \cite{fouloneng} for an English version.)\\\\
Antony Manning \cite{manning79} noticed that on non positively curved Riemannian manifolds, the topological entropy of the geodesic flow was equal to the volume entropy of the Riemannian metric. The volume entropy of a Riemannian metric $g$ on $M$ measures the exponential asymptotic growth of the volume of balls in the universal covering $\tilde M$ ; it is defined by 
$$h_{vol}(g)=\lim_{r\to\infty} \frac{1}{r} \log vol(B(x,r)),$$
where $vol$ denotes the Riemannian volume corresponding to $g$. We can also consider the volume entropy $h_{vol}(\o,d_{\o})$ of a Hilbert geometry $(\o,d_{\o})$ and extends the result of Manning in this case. This yields to the following rigidity result :

\begin{cor}\label{vol}
Let $\o$ be a strictly proper convex open set in $\P(\R^n)$ divided by a group $\Gamma\in PGL(\R^n)$ such that $M=\o\backslash \Gamma$ is compact. Then 
$$h_{vol}(\o,d_{\o})\leq n-1$$
with equality if and only if $\o$ is an ellipsoid.
\end{cor}
Thus, in the case of a manifold which admits a hyperbolic structure, the maximum of the (topological or volume) entropy characterizes the Teichm\"uller space $\mathcal{T}(M)$ in $\G(M)$. In any case, we get a function entropy $h : \G(M) \longrightarrow \R$ which takes its values in $(0,n-1]$. That yields to some natural questions :
\begin{itemize}
\item what is the infimum of $h$ and is it attained ?
\item in the case of a manifold which does not admit any hyperbolic structure, what is the supremum of $h$ and is it attained ?
\item how regular is $h$ ?\\
\end{itemize}

Let us now explain the contents of the paper.\\
After some necessary preliminaries consisting of basic facts, notations and motivations, we extend in section \ref{objects} the dynamical formalism introduced in \cite{foulon86} to our context. In particular, it allows us to define a notion of curvature, that we prove to be constant and strictly negative, and also to make parallel transport along the orbits of the geodesic flow.\\
In section \ref{sectionpaa}, this parallel transport is related to the action of the geodesic flow, that yields to a new description of the Anosov property which we guess is more intrinsic. Here the projective flatness of the structures is crucial : working in the universal covering identified with $\o$, we can indeed compare the parallel transport with respect to the Hilbert metric with the euclidean one (section \ref{sectionparom}) ; then an acute study allows to control the asymptotic behavior of the flow on the tangent space. This part is the technical core of the paper.\\
Using ergodic properties of hyperbolic systems and some arguments of symmetry, sections \ref{sectionlyap} and \ref{symmetry} prove the upper bound in theorem \ref{majeur}. Motivations and ideas of the proof appear in the preliminaries, section \ref{entropies}. These sections also give links between these dynamical properties, namely Lyapunov exponents, the group $\Gamma$ and the boundary of the convex $\o$.\\
Section \ref{equality} explicits the case of equality in theorem \ref{majeur} and provide some complementary facts and considerations about invariant measures. It also gives a large lower bound for the topological entropy in terms of regularity of the boundary of the convex.\\
Finally, the last section extends the results obtained by Manning, which leads to corollary \ref{vol}.\\

I would like to thank Patrick Foulon for all the interesting and fruitful discussions and ideas, Constantin Vernicos for his constructive remarks, Thomas Barthelm\'e and Camille Tardif for listening to my (sometimes strange) interrogations, and also Internet without which Ludovic Marquis, Yves Benoist, Jimmy Lamboley, Fran\c cois Ledrappier, Gerhard Knieper, and Fran\c cois Labourie could not have answered my questions.


\vskip 1cm

\section{Preliminaries : concepts and notations}

\subsection{Hilbert geometry}\label{hilbertgeom}
\subsubsection{Generalities}
\begin{figure}[h!]
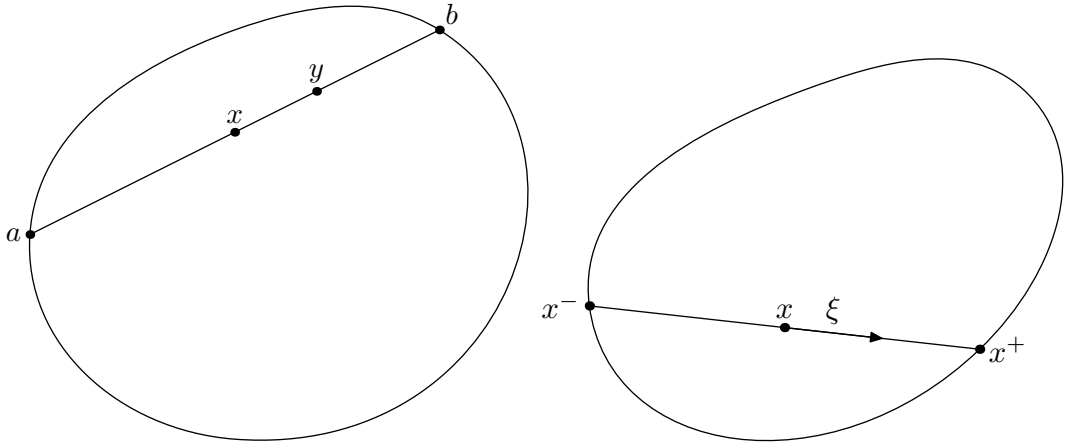

\begin{center}
\includegraphics[angle=0,width=7cm]{touteslesfigures.0} 
\includegraphics[angle=0,width=7cm]{touteslesfigures.1}

\end{center}
\caption{The Hilbert distance and the associated Finsler metric}
\label{fighilbert}
\end{figure}

Hilbert geometries were introduced (surprisingly) by David Hilbert as an example for what is now known as Hilbert's fourth problem : roughly speaking, characterize the metric geometries whose geodesics are straight lines. Hilbert geometries are defined in the following way.\\
Take a properly convex open set $\o$ of the projective space $\P^n(\R),\ n\geq 2$, where properly convex means you can find an affine chart in which $\o$ appears as a bounded convex set. The Hilbert metric $d_{\o}$ is defined in such an affine chart by 
$$d_{\Omega}(x,y)=\frac{1}{2}|\log([a,b,x,y])|,\ x,y \in \o,$$
where $a,b$ are the intersection points of the line $(xy)$ with the boundary $\partial \Omega$ ; $[a,b,x,y]=\frac{ax/bx}{ay/by}$ denotes the cross ratio of the four points (c.f. Figure \ref{fighilbert}).

Since cross-ratios are preserved by homographies, this definition does not depend on the choice of the affine chart. The space $(\o,d_{\o})$ is then a complete metric space ; see \cite{sm} for subsequent details.\\
In general, the metric is not Riemannian but Finslerian : instead of a quadratic form, we only have a convex norm on each tangent space. At the point $x\in \o$, this norm is given by
\begin{equation}\label{finsler}
F(x,\xi) = \frac{|\xi|}{2}\big(\frac{1}{\ xx^+}+\frac{1}{\ xx^-}\big),
\end{equation}
where $|\ .\ |$ denotes the euclidean norm, $x^+,\ x^-$ are the intersections of the line $\{x+\lambda \xi\}_{\lambda\in\R}$ with the boundary $\partial \Omega$ (see again Figure \ref{fighilbert}). From this formula, we see that $F:TM\longrightarrow M$ has the same regularity as the boundary $\partial\o$. 
Among all these geometries, those given by ellipsoids are particular : they are the only cases where the metric $F$ is Riemannian (see \cite{sm} for more precise statements), and in this case, $(\o,d_{\o})$ is nothing else than the Klein model for the hyperbolic space. Thus, a relevant problem is to compare the space $(\o, d_{\o})$ to standard geometries. In particular, note the two following opposite results.
\begin{thm}\label{hyp}
\begin{itemize}
\item \cite{cv} If $\o$ is $C^2$ with definite positive Hessian then the metric space $(\o,d_{\o})$ is bi-Lipschitz equivalent to the hyperbolic space $H^n$.
\item \cite{cvv} \cite{bernig} \cite{vernicos} $(\o,d_{\o})$ is bi-Lipschitz equivalent to the euclidean space if and only if $\o$ is a convex polytope, that is the convex hull of a finite number of points.
\end{itemize}
\end{thm}
From a different point of view, Benoist also found general conditions on the boundary for $(\o,d_{\o})$ to be Gromov-hyperbolic ; see \cite{benoistqs}. 

\subsubsection{Geodesics}
In any case, the space $(\o,d_{\o})$ is geodesically complete, where by geodesic, we mean a curve which locally minimizes the distance among all piecewise $C^1$ curves ; indeed,
any straight line is a geodesic. The converse is true if and only if the boundary of every plane section of $\o$ contains at most one open segment (see \cite{sm}). 

\subsubsection{Isometries}
The subgroup of elements of $PGL(n+1,\R)$ which preserve the convex $\o$ is obviously a subgroup of isometries of the space $(\o,d_{\o})$. The converse is false in general (see \cite{delaharpe}) and there is no known necessary and sufficient condition for this property to be true. The best sufficient condition was given in \cite{delaharpe} and specified in \cite{sm} : when the space is uniquely geodesic, then $Isom(\o,d_{\o}) \subset PGL(n+1,\R).$ In particular, this is true when $\o$ is strictly convex.

\subsection{Hilbert geometry on compact manifolds and divisible convex sets}

We now consider compact manifolds locally modeled on these geometries : we say that a manifold $M$ admits a convex projective structure if there exist a properly convex open set $\o$ and a subgroup $\Gamma\subset PGL(n+1,\R)$ preserving $\o$, such that $M=\o/\Gamma$. The convex set $\o$ is then said to be divisible. This structure identifies the universal covering of $M$ with $\o$, and its fundamental group $\pi_1(M)$ with $\Gamma$.\\
The ellipsoid is once more a particular case of divisible convex set. As was already noticed by Ehresmann, this is the only analytic model. In fact, for any divisible convex set which is not an ellipsoid, there exists some $0<\epsilon<1$ for which $\partial\o$ is not $C^{1+\epsilon}$. For more properties, especially about the groups $\Gamma$, look at the papers of Yves Benoist \cite{benoistcv1}, \cite{cv2}, \cite{cv3}, \cite{cv4}.\\
Among divisible convex sets, we have to distinguish the strictly convex and the non strictly convex ones ; indeed, if $\o$ is divisible by a group $\Gamma$ then the following are equivalents (\cite{benoistcv1}) :
\begin{itemize}
\item the space $(\o,d_{\o})$ is Gromov-hyperbolic ;
\item $\o$ is strictly convex ;
\item the boundary $\partial\o$ of $\o$ is $C^1$ ;
\item $\Gamma$ is Gromov-hyperbolic.
\end{itemize}
From these conditions, we see that all convex projective structures on a given manifold $M$ are either all strictly convex or all not strictly convex. In this paper, since we want to study the geodesic flow, we restrict ourselves to manifolds which admit strictly convex projective structures.\\
The set $\mathcal{G}(M)$ of projective convex structures is the set of equivalences classes of such representations of $M$, where two representations  $\o_1/\Gamma_1$ and   $\o_2/\Gamma_2$ are equivalent if there exists $g\in PGL(n+1,\R)$ such that $g(\o_1)=\o_2,\ g\Gamma_1 g^{-1} = \Gamma_2.$
Since in the strictly convex cases, $Isom(\o,d_{\o}) =\{g\in PGL(n+1,\R),\ g(\o)=\o\},$ looking at $\mathcal{G}(M)$ is equivalent as looking at the space of Hilbert metrics we can put on $M$ ; $\mathcal{G}(M)$ then plays the role that the Teichm\"uller space $\mathcal{T}(M)$ plays for hyperbolic metrics.\\
$\mathcal{T}(M)$ is naturally embedded in $\mathcal{G}(M)$. However, $\mathcal{G}(M)$ is really bigger than $\mathcal{T}(M)$, as already noticed in the introduction.

\subsection{Geodesic flow}\label{flow}
For every strictly convex projective structure on the compact manifold $M$, we are able to define the geodesic flow of the Hilbert metric since in this case, there is a unique geodesic between two points, which is a straight line in any projective chart.\\
In this paper we study the geodesic flow $\ph^t$ as defined on the homogeneous bundle 
$$HM = (TM\backslash \{0\}) / \R^*_+,$$
with projection $\pi : HM \longrightarrow M$ : a point $w=(x,[\xi]))\in HM$ is given by a point $x\in M$ and a direction $[\xi]$, where $\xi\in TM$. If $w=(x,[\xi]))\in HM$, then its image $\ph^t(w)=(x_t,[\xi_t]))$ is obtained by following the geodesic leaving $x$ in the direction $[\xi]$ during the time $t$, that is the length (for the Hilbert metric) of the corresponding geodesic curve between $x$ and $x_t$ is $t$ ; the direction $[\xi_t]$ is the direction tangent to this geodesic at the point $x_t$.\\
On the universal covering of $M$, identified with $\o$, the geodesic flow $\tilde{\ph}^t$ has a very simple interpretation : take a point $x\in \o$ and a direction $[\overrightarrow{xx^+}]$ for a point $x^+\in\partial\o$ ; the image $\tilde{\ph}^t(w)$ of $w=(x,[\overrightarrow{xx^+}])\in H\o$ by the geodesic flow is given by $(x_t,[\overrightarrow{x_tx^+}])$, where $d_{\o}(x,x_t)=t$. The flow $\ph^t$ on $HM$ is then obtained by using the projection $H\o \longrightarrow HM$.\\\\
The infinitesimal generator of the geodesic flow is a vector field $X$ defined on $HM$, that is a section $X : HM \longrightarrow THM$ of the tangent bundle of $HM$. On $H\o$, we thus get a $\Gamma$-invariant vector field $\tilde{X}$ ; since orbits of the flow are lines (the metric is said to be flat), there exists, once an affine chart is fixed, a function $m : HM \longrightarrow \R$ such that $\tilde{X}=mX^e,$ where $X^e$ denotes the infinitesimal generator of the euclidean metric on $\o\subset\R^n$. A direct calculation gives
$$m(x,[\xi]) = 2\left(\frac{1}{xx^+}+\frac{1}{xx^-}\right)^{-1} = 2\ \frac{xx^+\ xx^-}{x^+x^-}\ ,$$
so that $F(x,\xi)m(x,[\xi])=|\xi|.$ This property of flatness and the form of $m$ will be crucial to extend some concepts in section \ref{objects} despite the lack of regularity.\\

The geodesic flow of the Hilbert metric was studied by Yves Benoist, who proved the following
\begin{thmsansnum}[\cite{benoistcv1}]
The geodesic flow is a topologically mixing Anosov flow.
\end{thmsansnum}
Recall that a $C^1$ flow $\ph^t : W \longrightarrow W$ generated by $X$ on a compact manifold $W$ is an Anosov flow if there exist a decomposition
$$TW=\R.X \oplus E^s \oplus E^u,$$
and constants $C,\alpha,\beta>0$ such that for any $w\in W$ and $ t\geq 0$, 
$$\| d\ph^t(Z^s(w)) \| \leq C e^{-\alpha t},\ Z^s(w)\in E^s(w),$$
$$\| d\ph^{-t}(Z^u(w)) \| \leq C e^{-\beta t},\ Z^u(w) \in E^u(w).$$
Topologically mixing means that for any open sets $U,V \subset W$, there exists $T\geq 0$ such that for any $t\geq T$, $\ph^t(U)\cap V \not= \emptyset$.\\\\
Such a property was first established by Hadamard \cite{hadamard} in 1898 for the geodesic flow on hyperbolic surfaces, and then generalized to Riemannian manifolds of negative curvature by Anosov in the famous \cite{anosov}. It is thus a property that is shared by our geometries. Our goal is to study what dynamically separates Riemannian hyperbolic structures from the others; that is to find dynamical properties which characterize hyperbolic metrics among the non Riemannian Hilbert metrics.\\
Benoist made a first step by proving the
\begin{prop}[\cite{benoistcv1}, Proposition 6.7]\label{absco}
There exists an absolutely continuous $\ph^t$-invariant measure if and only if the Hilbert metric is Riemannian.
\end{prop}
Recall that a measure $\mu$ on a manifold $W$ is said to be absolutely continuous (or smooth) if it is in the Lebesgue class : if $A$ is a Borel subset of $W$, then $\mu(A)=0$ as soon as $\lambda(A)=0$, where $\lambda$ denotes a Lebesgue measure on $HM$.
The proposition above will be useful in section \ref{equality} to determine the case of equality in theorem \ref{majeur}.

\subsection{Topological and measure theoretic entropies}\label{entropies}

Let $\ph^t:W\longrightarrow W$ be a flow on a compact manifold $W$. For $t\geq 0$, we define the distance $d_t$ on $W$ by :
$$d_t(x,y)=\max_{0\leq s\leq t}\ d(\varphi^s(x),\varphi^s(y)),\ x,y\in W.$$
For any $\epsilon >0$ and $t\in \R$, we consider coverings of $W$ by open sets of diameter less than $\epsilon$ for the metric $d_t$. Let $N(\varphi,t,\epsilon)$ be the minimal cardinality of such a covering. The topological entropy (\cite{adler}) of the flow is then the  (well defined) quantity
$$h_{top}(\varphi)=\lim_{\epsilon \to 0}\Big[ \limsup_{t \to \infty}\frac{1}{t}\log N(\varphi,t,\epsilon)\Big].$$
In a certain sense, it measures how much the system is chaotic. In particular, the topological entropy of an Anosov flow is nonnegative. It appears in various and numerous contexts ; the most celebrated result may be this one, essentially due to Margulis (see \cite{margulisbook}, \cite{kh}) : if $\ph$ is a topologically mixing Anosov flow, then the number $P_T(\ph)$ of closed orbits of length less than $T$ satisfies the following asymptotic equivalent, with $h=h_{top}(\ph)$ :  
$$P_T(\ph) \sim \frac{e^{-hT}}{hT}.$$
As an example, the topological entropy of the geodesic flow of a compact hyperbolic manifold of dimension $n\geq 2$ is $(n-1)$. Our main theorem \ref{majeur} states that this property characterizes the hyperbolic structures among all strictly convex projective ones.\\\\
To prove this theorem, we will make use of certain objects and results that appear in the ergodic theory of hyperbolic dynamical systems. Here come the motivations for the proof.\\
Let $\M$ denote the set of $\ph^t$-invariant probability measures. To any $\mu\in \M$ is attached a number $h_{\mu}$ called measure-theoretic entropy ; for definition and basic properties, see \cite{kh} or \cite{walters}. The variational principle (\cite{goodman} or \cite{misiurewicz}) states that 
$$h_{top}(\ph)=\sup_{\mu\in\M} h_{\mu},$$
and in the case of a topologically mixing $C^{1+\epsilon}$ Anosov flow (that is relevant for us), we know from Bowen \cite{bowen1} and/or Margulis \cite{margulis} (see also \cite{kh}) that there exists a unique measure $\mu_{BM}$, now known as the Bowen-Margulis measure, such that
$$h_{\mu_{BM}}=h_{top}(\ph).$$
On a hyperbolic manifold, the Bowen-Margulis measure of the geodesic flow is the natural Liouville measure. From proposition \ref{absco}, we know that in the case of a non Riemannian Hilbert metric, this measure will not be smooth anymore.\\
Osedelec's theorem \cite{osedelec} and Pesin-Ruelle inequality \cite{ruelle78} give a way to  calculate $h_{\mu_{BM}}$ : if $\mu\in\M$ then the set of regular points is of full measure (see definition \ref{defilyapunov} and theorem \ref{fullmeasure}) and
\begin{equation}\label{ruelle}
h_{\mu}\leq \int \chi^+ d\mu,
\end{equation}
where $\chi^+$ is the sum of positive Lyapunov exponents. Proposition \ref{core} will give a formula for our Lyapunov exponents which appears to be sufficient to conclude.

\subsection{Volume entropy of Hilbert geometries}

We define the volume entropy of a Hilbert geometry $(\o,d_{\o})$, provided it exists, by
\begin{equation}\label{entvol}
h_{vol}(\o,d_{\o})=\lim_{r\to\infty} \frac{1}{r} \log vol(B(x,r)).
\end{equation}
It measures the asymptotical exponential growth of the volume of balls. By volume, we mean the Hausdorff measure associated to the Hilbert metric. Note that, if the convex set is divisible by a group $\Gamma$, this volume is $\Gamma$-invariant, giving a volume on the manifold $\o/\Gamma$.\\
The problem of measuring a volume in a Finsler space was already discussed a lot and we will not discuss it again. Look at \cite{bubu} and \cite{alvarez} for instance.\\\\
It is not clear when the limit in (\ref{entvol}) exists, but some results are already known : as a consequence of theorem \ref{hyp}, if $\o$ is a polytope then $h_{vol}(\o,d_{\o})=0$ ; at the opposite, we have the
\begin{thm}[\cite{bbv}]
If $\o$ is $C^{1,1}$ then $h_{vol}(\o,d_{\o})=n-1$. 
\end{thm} 
It is conjectured that $h_{vol}(\o,d_{\o})\leq n-1$ for any convex set $\o$ of dimension $n$. In \cite{bbv} was proved the conjecture in dimension $n=2$ and was also constructed an example with $0<h_{vol}<1$. Theorem \ref{vol} will provide numerous examples of convex sets, in any dimension $n$, whose entropy satisfies 
$$0<h_{vol}<n-1.$$


\vskip 1cm 
\section{Dynamical formalism}\label{objects}

To prove the main theorem, we use the dynamical objects introduced by Patrick Foulon in \cite{foulon86} to study second order differential equations : they provide analogues of Riemannian objects such as covariant differentiation, parallel transport and curvature for any such equation which is regular enough. Here we deal with a more irregular case but the objects are still well defined, and even smooth, when we restrict ourselves to move along an orbit.\\
In what follows, $M$ is a compact manifold with a strictly convex projective structure given by $M=\o/\Gamma$. The notations are those of the preliminaries.

\subsection{Directional smoothness}

Assume a complete vector field $X$ is given on a manifold $W$. We denote by $C^{p}_X(W)$ (or simply $C^{p}_X$) the set of $C^p$ functions $f$ on $W$ which are smooth in the direction $W$, that is $L^n_Z f$ exists and is continuous for any $n\geq 1$.\\
We say that a $C^p$ vector field $Z$, $p\geq 0$, is smooth in the direction $X$, or along any orbit of the flow of $X$, if $Z$ can be locally written as $Z=\sum f_i Z_i$ where the $Z_i$ are smooth vector fields on $W$, and $f_i\in C^p_X$.

\subsection{Some objects}

Here we introduce a continuous decomposition of the tangent space $THM$, which is smooth along any orbit of the flow of the Hilbert metric, generated by $X$. We do it on $H\o$ and then come back to $M$ by using the projection $H\o\longrightarrow HM$. On $H\o$, we have two vector fields, $\tilde X$ and $X^e$ related by $\tilde X=mX^e$.\\
The key remark is that the function $m$ on $HM$ is smooth in the direction $\tilde X$ (or, equivalently, $X^e$). More precisely, we have
for $w=(x,[\xi])\in H\o$, 
$$L_{X^e} m\ (w) = 2\ \frac{xx^+ - xx^-}{x^+x^-}\ ;\ L_{X^e}^2 m\ (w) = -\frac{4}{x^+x^-},\ L^n_{X^e} m = 0,\ n\geq 3.$$
From this, we also see that $L_{\tilde X} m$ itself is $C^1$. Thus, for any vector field $Z$ on $H\o$ such that $X(w)$ and $Z(w)$ are nowhere collinear, $w\in H\o$, $L_Z (L_X m)$ exists and is continuous ; hence $L_X (L_Z m)$ also from the following version of Schwarz' theorem.

\begin{lemma}
Let $f : \R^2 \longrightarrow \R$ be a $C^1$ map. If $\frac{\partial^2 f}{\partial x\partial y}$ exists and is continuous then so is  $\frac{\partial^2 f}{\partial y\partial x}$. 
\end{lemma}

\subsubsection{Vertical vectors and the verticality lemma} 
The vertical distribution is the smooth distribution $VH\o = \ker d\pi$ where $\pi : H\o \longrightarrow \o$ is the bundle projection. The letter $Y$ will always denote a vertical vector field along $\tilde\ph.w$. The following lemma is proved in \cite{foulon86} :

\begin{lemma}
Let $Y_1,\cdots,Y_{n-1}$ be a base of $VH\o$ along $\tilde\ph.w$. Then the family
$$\tilde X,Y_1,\cdots,Y_{n-1},[\tilde X,Y_1],\cdots,[\tilde X,Y_{n-1}]$$
is a base of $TH\o$.
\end{lemma}

\subsubsection{The vertical endomorphism}
From this lemma, we can define on each tangent space $T_p H\o$ the vertical endomorphism $v_{\tilde X}(p)$ given by :
\begin{itemize}
\item $v_{\tilde X}(p)(\tilde X)= v_{\tilde X}(p)(Y)=0;$
\item $v_{\tilde X}(p)([\tilde X,Y])=-Y(p)$,
\end{itemize}
which is a kind of projection on the vertical subspace $V_p H\o$. This allows to define a vertical operator $v_{\tilde X}$ on $TH\o$ by setting $v_{\tilde X}(Z)(p)=v_{\tilde X}(p)(Z)$, such that for any function $f$ on $H\o$, we have $v_{\tilde X}(fZ)=fv_{\tilde X}(Z)$. From the very definition, we check that $v_{\tilde X}=m v_{X^e}$.\\

\subsubsection{Horizontal considerations}
The horizontal operator $H_{\tilde X} : VH\o \longrightarrow T\ho$ is then defined by :
$$H_{\tilde X}(Y) = -[\tilde X,Y] - \frac{1}{2} v_{\tilde X} ([\tilde X,[\tilde X,Y]]).$$
From the introductory remarks about $m$, $L_{\tilde X}L_Y m$ exists and is continuous, so the Lie brackets
$$[\tilde X,Y]=m[X^e,Y] - L_Y m X^e$$
and
$$[\tilde X,[\tilde X,Y]] = m^2 [X^e,[X^e,Y]] + L_{\tilde X} m [ X^e,Y] - (L_{\tilde X}L_Y m -m L_{[\tilde X,Y]}m )  X^e.$$
are well defined ; we finally get, since $v_{\tilde X}=m v_{X^e}$,
\begin{equation}\label{horizontal} H_{\tilde X}(Y) = mH_{X^e}(Y) + L_Y(\log m) \tilde{X} + \frac{1}{2}L_{\tilde X}(\log m) Y. \end{equation}
The operator $H_{\tilde X}$ is linear : for any function $f\in C^0_{\tilde X}$, we can compute $H_{\tilde X}(fY)$ and we have
$$H_{\tilde X}(fY)= fH_{\tilde X}(Y).$$
The horizontal distribution $h^{\tilde X}H\o$ is then defined as the image of $VH\o$ by $H^{\tilde X}$. By the help of the verticality lemma, we can prove that $H_{\tilde X}$ is injective and we get the continuous decomposition  
$$TH\o=\R.\tilde X \oplus VH\o \oplus h^{\tilde X}H\o$$
which, as shown by $(\ref{horizontal})$, is smooth along any orbit, that is in the direction $\tilde X$. The operators $v_{\tilde X}$ and $H_{\tilde X}$ exchange $VH\o$ and $h^{\tilde X}H\o$ : the constructions above allow us to consider the compositions 
$v_{\tilde X} \circ H_{\tilde X}$ and $H_{\tilde X} \circ v_{\tilde X}$, and see that 
$$v_{\tilde X} \circ H_{\tilde X} = Id_{VH\o}, \ (H_{\tilde X} \circ v_{\tilde X})_{h^{\tilde X}H\o} = Id_{h^{\tilde X}H\o}.$$ 
Setting $J^{\tilde X}=v_{\tilde X}$ on $h^{\tilde X}H\o$ and $J^{\tilde X}=-H_{\tilde X}$ on $VH\o$ defines a complex structure on $h^{\tilde X}H\o \oplus VH\o$.\\

\paragraph{Notations :}
\begin{itemize}
\item We associate to the decomposition
$$TH\o=\R.\tilde X \oplus VH\o \oplus h^{\tilde X}H\o$$
the corresponding decomposition of the identity :
$$Id = p^{\tilde X} \oplus p_v^{\tilde X} \oplus p_h^{\tilde X}.$$
\item By a horizontal vector field, we will mean a vector field $h\in h^{\tilde X}H\o$, such that $h=H_{\tilde X}(Y)$ for a certain vertical vector field $Y$ ; any such $h$ is smooth along $\tilde X$.
\end{itemize}
Let us remark that $p_h^{\tilde X} = H_{\tilde X}\circ v_{\tilde X}$. Moreover,
\begin{lemma}\label{projection}
 For any smooth vector field $Z$, we have
$$p_{\tilde X}(Z)=p_{X^e}(Z) -  L_{v_{X^e}(Z)}(\log m) X^e;$$
$$p_v^{\tilde X}(Z)= p_v^{X^e}(Z)-\frac{1}{2} (L_{X^e} \log m) v_{X^e}(Z); $$
$$p_h^{\tilde X}(Z)=p_h^{X^e}(Z) + (L_{v_{X^e}(Z)} (\ln m)) X^e + 
\frac{1}{2} (L_{X^e} \log m) v_{X^e}(Z).$$
In particular, every projection of $Z$ is still smooth along $\tilde X$.
\end{lemma}
\begin{proof}
Let $Z=a{\tilde X}+Y+h=a^eX^e+Y^e+h^e$ be the two decompositions of the vector field $Z$ along $\tilde\ph.w$.
If we note $y=v_{X^e}(h^e)=v_{X^e}(Z)$, we have by using (\ref{horizontal})
$$h=H_{\tilde X}(v_{\tilde X}(Z))= \frac{1}{m} H_{\tilde X}(y)= H_{X^e}(y) + \frac{1}{2m}L_{\tilde X}(\log m) y + \frac{1}{m}L_y(\log m) \tilde{X}.$$
Thus $$h=h^e + \frac{1}{2}L_{X^e}(\log m) y + L_y(\log m) X^e,$$
and $$Z=(a \tilde X+ L_y(\log m) X^e)+(Y+\frac{1}{2}L_{X^e}(\log m) y)+h^e=a^eX^e+Y^e+h^e.$$
Identifying gives the result.
\end{proof}

\subsubsection{Dynamical derivation}
We define an analog of the covariant derivation along $\tilde X$ that we call the dynamical derivation and note $D^{\tilde X}$. It is a differential operator of order 1 : if $f\in C^0_{\tilde X}$, then $D^{\tilde X}(fZ)=fD^{\tilde X}(Z) + (L_{\tilde X}f) Z$. $D^{\tilde X}$ can be defined by setting 
$$D^{\tilde X}(\tilde X)=0, \ D^{\tilde X}(Y)=\displaystyle -\frac{1}{2} v_{\tilde X}([\tilde X,[\tilde X,Y]]), \ [D^{\tilde X} , H_{\tilde X}]= 0.$$
Note that $[D^{\tilde X} , v_{\tilde X}]= 0$ ; that, on $VH\o$, we can write 
$$D^{\tilde X}(Y) = H_{\tilde X}(Y) +[\tilde X,Y];$$
and that we have 
$$D^{\tilde X}=mD^{X^e} + \frac{1}{2} L_{\tilde X}(\log m) Id.$$
Also observe that $D^{X^e}(Z)$ is the usual covariant derivative $\nabla_{X^e} Z$ of a vector field $Z$ along $X^e$.\\
A vector field $Z$ is said to be parallel along $\tilde X$, or along any orbit, if $D^{\tilde X} (Z)=0$. In section \ref{sectionpaa}, we will consider parallel transport with respect to this dynamical derivation. In the euclidean case, the projection on the base of this transport coincides with the usual parallel transport along geodesics.

\subsubsection{Jacobi endomorphism and curvature}
The Jacobi operator $R^{\tilde X}$ associated to $\tilde X$ on $TH\o$ is a linear operator which respect to any function in $C^0_{\tilde X}$ ; it is defined by 
$$R^{\tilde X}(\tilde X)=0, \ R^{\tilde X}(Y)= p_v^{\tilde X}([\tilde X,H_{\tilde X}(Y)]), \ [R^{\tilde X} , H_{\tilde X}]= 0,$$
which is well defined thanks to lemma \ref{projection}.
Note that $[R^{\tilde X} , v_{\tilde X}]= 0$ and on $VH\o$ :
$$R^{\tilde X}=m^2R^{X^e} + \big(\frac{1}{2} L^2_{\tilde X}(\log m) -\frac{1}{4}(L_{\tilde X}(\log m))^2 \big) Id.$$
This is an analogous of Riemannian curvature in our context. In particular, a direct calculation gives
$R^{X^e}=0$ and $R^{\tilde X}|_{VH\o}=-Id|_{VH\o}$.

\subsubsection{On the manifold}
All these objects with respect to $\tilde X$ can now be transported on the homogeneous tangent space $HM$ of the given manifold $M$. We thus have a continuous decomposition
$$THM = \R.X \oplus VHM \oplus h^XHM$$
which is smooth along any orbit of $X$, and a complex structure given by $J^X$ on $VHM \oplus h^XHM$. The dynamical derivation $D^X$ will be the central object in the two next sections. The restriction of the Jacobi endomorphism to $VHM\oplus h^XHM$ is still $R^X|_{VHM\oplus h^XHM}=-Id$ : the Finsler manifold $(M,F)$ has then constant strictly negative curvature.

\subsection{The 1-form associated to a Finsler metric}
The vertical derivative of a Finsler metric $F$ on a manifold $M$ is the 1-form on $TM$ defined for $Z\in T(TM)$ by : 
$$d_vF(x,\xi)(Z) = \lim_{\epsilon\rightarrow 0} \frac{F(x,\xi+\epsilon dp(Z))-F(x,\xi)}{\epsilon},$$
 where $p : TM \longrightarrow M$ is the projection on the base. This form depends only on the direction $[\xi]$, that is $d_vF(x,\lambda\xi)(Z)=d_vF(x,\xi)(Z)$ for $\lambda>0$. As a consequence, it defines a 1-form $A$ on $HM$ by using the projection 
$TM\backslash \{0\} \longrightarrow HM$. Let $X$ be the infinitesimal generator of the geodesic flow of $F$ on $HM$. Since $[d\pi(X(x,[\xi]))]=[\xi]$, we can define $A$ for any $Z\in THM$ by 
$$A(Z) = \lim_{\epsilon\rightarrow 0} \frac{F(d\pi(X+\epsilon Z)) -1}{\epsilon}.$$
Note that $A(X)=1$ and that $A(Y)=0$ for any vertical vector field. The following proposition is well known for regular Finsler metrics and in this case easier to prove since we are allowed to derivate $A$ to get the 2-form $dA$.\\\\

\begin{prop}\label{kerA} Let $M=\o/\Gamma$ be a compact manifold with a strictly convex projective structure and $A$ the 1-form on $HM$ associated to the Hilbert metric on $M$. Then
$$\ker A= VHM \oplus h^XHM.$$
Furthermore, $A$ is invariant under the geodesic flow of the Hilbert metric.
\end{prop}

We prove the proposition on $H\o$ to where we can do some calculus, and by use of the covering map $\Omega \rightarrow M$ we get the result on $HM$. Choose a point $w=(x,[\xi])\in H\o$ with orbit $\tilde\ph .w$. We will work in {\it an affine chart adapted for this orbit}, where the intersection space $T_{x^+}\partial\o\cap T_{x^-}\partial\o$ is at infinity, so that $T_{x^+}\partial\o$ and $T_{x^-}\partial\o$ are parallel, and orthogonal to $\overrightarrow{xx^+}$. 
\begin{center} {\it All along this paper, when we talk about a {\bf good affine chart} or a {\bf chart adapted} at $w \in H\o$ or its orbit $\tilde\ph.w$, we mean such an affine chart. (See Figure \ref{figgoodchart})}
\begin{figure}[h!]
\begin{center}
\includegraphics[angle=0,width=9cm]{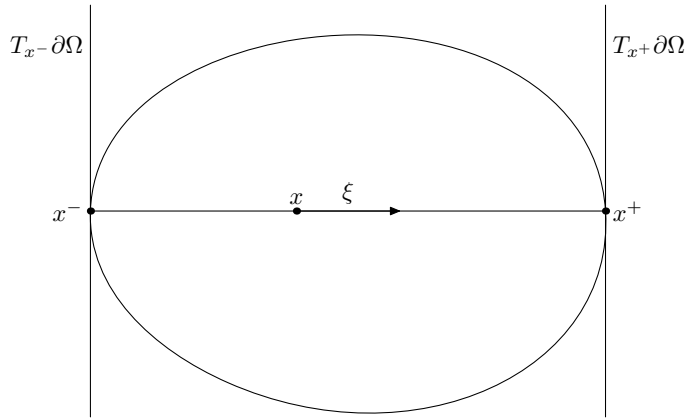}
\end{center}
\caption{A good chart at $w=(x,[\xi])$}
\label{figgoodchart}
\end{figure}
\end{center}
In a good affine chart at $w$, we clearly have $L_Y m =0$ along $\tilde\ph.w$ for any vertical vector field $Y$. The proof will in fact show that $\ker dm = VH\o \oplus h^{\tilde X} H\o$ along $\tilde\ph.w$. In particular, we will have 
$$d\pi(V_wH\o\oplus h_w^{\tilde X}H\o) = \left(\overrightarrow{xx^+}\right)^{\perp}.$$

\begin{proof}[Proof of proposition \ref{kerA}] We only have to prove that $h^XH\o\subset \ker A$. To do some explicit differential calculus on $H\o$, we need to choose good coordinates on $H\o$. As we deal with directions and lines, instead of using an identification of $H\o$ with a unitary tangent bundle $T^1\o$, we prefer to use a kind of projective charts.\\
Let $w_0=(x_0,[\xi_0])$ be any point in $H\o$ and choose a small open neighborhood $U$ of $w_0$ in $H\o$. Assume the affine chart for $\o$ is adapted to $w_0$. If the neighborhood is small enough, we can choose the following coordinates : 
\begin{itemize}
\item $w_0=(x_0,[\xi_0])$ is the origin : $x_0 = 0 \in \R^n$ and $[\xi_0]=[1:0:\cdots:0], \ \xi_0=\frac{\partial}{\partial x_1}\in S^{n-1}$, where we make use of homogeneous coordinates and $S^{n-1}$ is the euclidean sphere ;
\item for $w=(x,[\xi])\in U$, the coordinates of $x$ are the usual euclidean coordinates in $\R^n$, and every $[\xi]$ is written as $[1:\xi_2:\cdots:\xi_n]$, where the $\xi_i$ vary in a neighborhood of $0$.
\end{itemize}
We use the associated base $\left(\frac{\partial}{\partial x_i}, \frac{\partial}{\partial \xi_j}\right)_{1 \leq i \leq n, 2 \leq j \leq n}$ on the tangent space $TU\subset TH\o$. By $\xi \in T\o$ we denote the vector $$\xi=\frac{\partial}{\partial x_1} +  \sum_{i=2}^{n} \xi_i \frac{\partial}{\partial x_i}.$$
In this chart, 
\begin{itemize} 
\item we define the vector field $X^0$ on $H\o$ by
$$X^0(w)=X^0(x,[\xi])=\frac{\partial}{\partial x_1} +  \sum_{i=2}^{n} \xi_i \frac{\partial}{\partial x_i}.$$
In particular, we have $X^0(w_0)=\frac{\partial}{\partial x_1}$  and $d\pi(X^0(x,[\xi]))=\xi$; moreover $\tilde X$ can be written as
$\tilde X=m^{-1} X^0$ with the function $m$ defined on $TU$ by $m(w)=F(d\pi(X^0(w)))=F(x,\xi);$
\item the vertical distribution is here given by $$VU=\textrm{vect}\{\frac{\partial}{\partial {\xi_i}}\}_{i\in\{2,\cdots,n\}};$$
\item since $L_Y m=0$ along $\tilde\ph.w_0$, the complex structure along $\tilde\ph.w_0$ given by $X^0$ is very simple. We have$$[X^0,\frac{\partial}{\partial \xi_j}]=-\frac{\partial}{\partial x_j},\ \ 
[X^0,[X^0,\frac{\partial}{\partial \xi_j}]]=0,\ j=2,\cdots,n,$$
hence
$$v_{X^0}(\frac{\partial}{\partial x_j})=\frac{\partial}{\partial \xi_j},\ \ 
H_{X^0}(\frac{\partial}{\partial \xi_j})=\frac{\partial}{\partial x_j},\ j=2,\cdots,n,$$
thus
$$h^{X^0}HU=\textrm{vect}\{\frac{\partial}{\partial {x_i}}\}_{i\in\{2,\cdots,n\}}.$$
\end{itemize}
Then, from (\ref{horizontal}), any horizontal vector field $h\in h^{\tilde X}U$ along $\tilde\ph.w_0$ can be written 
$$h=m^{-1} H_{X^0}(Y)- \frac{1}{2} (L_{\tilde X} \log m) Y,$$
for a certain vector $Y\in VU$. Note $h^0=H_{X^0}(Y)$. Then $A(h)=A(h_0)$, so we only have to prove that for any $i\in\{2,\cdots,n\}$, $A(\frac{\partial}{\partial {x_i}})=0$. But
$$A(\frac{\partial}{\partial x_i})= 
\lim_{\epsilon\rightarrow 0} \frac{F(d\pi(X+\epsilon \frac{\partial}{\partial x_i})) -1}{\epsilon}=
\lim_{\epsilon\rightarrow 0} \frac{F(d\pi(X^0+\epsilon \frac{\partial}{\partial {x_i}})) -F(d\pi(X^0))}{\epsilon}$$
so that, for $w\in\tilde\ph.w_0$,
$$A(w)(\frac{\partial}{\partial x_i})=
\lim_{\epsilon\rightarrow 0} \frac{F(x,\xi+\epsilon \frac{\partial}{\partial x_i})) -F(x,\xi)}{\epsilon} = 
D_{(x,\xi)}F(\frac{\partial}{\partial x_i}),$$
where we see $F$ as a real valued function on $\o \times \R^n \subset \R^{2n}$ with coordinates $(x_1,\cdots,x_n,\frac{\partial}{\partial x_1},\cdots,\frac{\partial}{\partial x_n})$. But in this chart, from the formula giving $F$, we clearly have $\frac{\partial}{\partial x_i}\in\ker DF,\ i=2 \cdots n,$ which proves that 
$h^{\tilde X}H\o\subset \ker A$. Hence, coming back to $M$, we get $h^XHM\subset \ker A$.\\

To prove that $A$ is invariant under the flow, we only have to prove that its kernel is invariant, which from the first result is equivalent to proving that
$$p_X([X,Y])=p_X([X,h])=0$$
for any vertical and horizontal vector fields $Y$ and $h$. Since $[X,Y]=-H_X(Y) + D^X(Y)$, we clearly have $p_X([X,Y])=0$. Furthermore, once again working in $H\o$ in a chart adapted to an orbit $\tilde\ph.w_0$, we have from lemma \ref{projection}, $p_{\tilde X}=p_{X^0}$; hence
$$p_{\tilde X}([\tilde X,h]) = p_{X^0} (m^{-1}[X^0,h] - L_h m^{-1} X^0)= m^{-1}p_{X^0} ([X^0,h])  - L_h m^{-1}.$$
But, in this affine chart, we also have $L_h m^{-1} =0$ : this can be seen directly or using the coordinates that were introduced before. Moreover, if $h=H_{\tilde X}(Y)$ and $h^0=H_{X^0}(Y)$, then
$$p_{X^0} ([X^0,h]) = p_{X^0} ([X^0,m^{-1} h^0 - \frac{1}{2} (L_{\tilde X} \log m) Y]) =
m^{-1} p_{X^0} ([X^0, h^0]) = 0.$$
Finally $p_{\tilde X}([\tilde X,h])=0$.
\end{proof}


\vskip 1cm
\section{Parallel transport and the Anosov property}\label{sectionpaa}

\subsection{Action of the flow on the tangent space}
We pick a tangent vector $Z(w)\in T_{w}HM$. We want to study the behavior of the vector field $Z(\ph^t(w))=d\ph^t(Z(w))$ defined along the orbit $\ph.w$. Assume $$Z(w)=Y(w)+h(w)\in V_{w}HM\oplus h^X_{w}HM.$$ Since $VHM\oplus h^XHM$ is invariant under the flow, we can write $Z=Y+h$. To find the expressions of $Y$ and $h$, we write that, since $Z$ is invariant under the flow, the Lie bracket $[X,Z]$ is 0 everywhere on $\ph.w$.\\
For that, let $(h_1,\cdots,h_{n-1})$ be a base of $h^XHM$ of $D^X$-parallel vectors along $\ph.w$, that is $h^t_i=h_i(\ph^t(w))=T^t(h_i(w))$, where $T^t$ denotes the parallel transport for $D^X$, $(h_i(w))_i$ is a base of $h_w^XHM$. Since $D^X$ and $v_X$ commute, the family $\{Y_i\}=\{v_X(h_i)\}$ is a base of $VHM$  of $D^X$-parallel vectors along $\ph.w$. We have immediately $h_i=H_X(Y_i)$ and we can check that 
\begin{equation}\label{bracket} [X,h_i]=-Y_i;\ [X,Y_i]=-h_i. \end{equation}
Then, in this base, $Z$ can be written as
$$Z=\sum a_i h_i + b_i Y_i,$$
where $a_i$ and $b_i$ are smooth real functions along $\ph.w$. The formulas (\ref{bracket}) give
$$\begin{array}{llll}
[X,Z]=0 & \Longleftrightarrow & \sum (L_X a_i - b_i) h_i + (L_X b_i-a_i)Y_i \\
& \Longleftrightarrow & b_i=L_X a_i ;\ a_i=L_X b_i,\ i=1, \cdots, n-1\\
& \Longleftrightarrow & b_i=L_X a_i ;\ a_i-L^2_X a_i,\ i=1,\cdots,n-1.
\end{array}$$
From that we get the solution
\begin{equation}\label{flot}
Z(\ph^t(w))=d\ph^t(Z(w))=\sum A_i e^t(h^t_i+Y^t_i) + B_i e^{-t} (h^t_i-Y^t_i),
\end{equation}
where
$$A_i=\frac{1}{2}(a_i(w)+b_i(w)),\ B_i=\frac{1}{2}(a_i(w)-b_i(w))$$
depend on the initial coordinates of $Z$ at $w$.\\\\

\subsection{The Anosov property}

Let us define the two diagonals $E^u$ and $E^s$ by
$$E^u=\{Y+H_X(Y), Y\in VHM\}, E^s=\{Y-H_X(Y), Y\in VHM\}=J^X(E^u).$$
We see from (\ref{flot}) that $E^u$ and $E^s$ are invariant under the flow. Furthermore if $Z^s(w)\in D^s(w),\ Z^u(w)\in D^u(w)$, then 
\begin{equation}\label{dd} d\ph^t(Z^u(w))= e^{t} T^t(Z^u(w)),\ d\ph^t(Z^s(w))= e^{-t} T^t(Z^s(w)). \end{equation}

\begin{thm}\label{anosov}
The geodesic flow $\ph^t$ is an Anosov flow with decomposition
$$THM=\R.X \oplus E^s \oplus E^u,$$
that is there exist constants $C,\alpha,\beta>0$ such that for any $w\in HM$ and $ t\geq 0$, 
$$\| d\ph^t(Z^s(w)) \| \leq C e^{-\alpha t},\ Z^s(w)\in E^s(w),$$
$$\| d\ph^{-t}(Z^u(w)) \| \leq C e^{-\beta t},\ Z^u(w) \in E^u(w).$$
\end{thm}
To prove this theorem, equations (\ref{dd}) above motivate the study of the parallel transport $T^t$ along an orbit  : we will thus focus on the exponential behavior of 
$\|T^t(Z^u(w))\|$ and $\|T^t(Z^s(w))\|$. The proof of the theorem will be completed in section \ref{proofanosov}.



\subsection{Comparison lemma}
Here is the key lemma, essentially due to Yves Benoist \cite{benoistcv1}. We note $E^{u,s}=E^u\cup E^s$.
\begin{lemma}\label{comparison}
For any Riemannian metric $\|. \|$ on $HM$, there exists a constant $C>0$ such that for any 
$Z(w)\in  E^{u,s}(w)$, 
$$C^{-1}\| Z(w)\| \leq F(d\pi(Z(w))) \leq C\| Z(w)\|.$$
\end{lemma}
\begin{proof}
Since $F:TM\rightarrow [0,+\infty)$ is a continuous function, so is the function
$$\begin{array}{rccc}
F\circ d\pi : &(E^s,\|.\|) &\longrightarrow & [0,+\infty[ \\
&u & \longmapsto & F\circ d\pi(u)
\end{array}$$ 
Thus its restriction to the compact $E^s_1=\{u\in E^s,\ \|u\|=1\}$ is bounded. Since it is also non zero, there exists $C>0$ such that, for any $u\in E^s_1$,
$$\frac{1}{C} \leq F(x,d\pi(u)) \leq C,$$
and we conclude the proof using the homogeneity of $F$. The same works for $E^u$.
\end{proof}

This lemma gives a way to tackle the problem : for any $Z(w)=Y(w)+h(w) \in E^{u,s}(w)$, the exponential behavior of $\|T^t(Z(w))\|$ will be the same as the one of $F(d\pi(T^t(h(w)))).$ 

\subsection{Parallel transports on $H\o$}\label{sectionparom}

We now come back on $\ho$ where we can do some calculus. The Riemannian metric $\|. \|$ and the Finsler metric $F$ on $HM$ give $\Gamma$-invariant metrics on $\ho$, that we also write $\|. \|$ and $F$. The lemma \ref{comparison} is still valid.\\
On $\ho$ we work with two vector fields, namely $\tilde X$ and $X^e$, with $\tilde X=mX^e$. $\tilde T^t$ and $T_e^t$ will denote respectively $D^{\tilde X}$ and $D^{X^e}$ parallel transports ; $\tilde E^s,\ \tilde E^u$ and $\tilde E^{u,s} \subset TH\o$ correspond to $E^s,\ E^u$ and $E^{u,s}$.

\begin{lemma}\label{horver}
If $Y(w)\in V_wH\o$ then 
$$\tilde T^t(Y(w))= \left(\frac{m(w)}{m(\ph^t(w))}\right)^{1/2} T_e^t(Y(w)).$$
Furthermore, in a good affine chart at $w$, if $h(w)\in h_w^{\tilde X}\ho$ then
$$d\pi(\tilde T^t(h(w)))=-(m(w)m(\ph^t(w)))^{1/2} d\pi(T_e^t(h(w))).$$
\end{lemma}
\begin{proof}
 We look for the unique vector field $Y$ along $\tilde \ph.w$ such that $D^{\tilde X}(Y)=0$ and which is  equal to $Y(w)$ at the point $w$. We recall that
$$D^{\tilde X}(Y)=mD^{X^e}(Y) + \frac{1}{2} L_{\tilde X}(\log m) Y.$$
Assume we can write $Y=fY^e$, where $Y^e$ is parallel for $D^{X^e}$ along $\tilde \ph.w$. Then $f$ is the solution of the equation
$$L_{\tilde X} (\log f) + \frac{1}{2} L_{\tilde X}(\log m)=0,$$
which with $f(w)=1$ gives
$$f(\tilde \ph^t(w))=\left(\frac{m(w)}{m(\tilde \ph^t(w))}\right)^{1/2}.$$
Finally,
$$\tilde T^t(Y(w))=\left(\frac{m(w)}{m(\tilde \ph^t(w))}\right)^{1/2} T_e^t(Y(w)).$$

Now, let $h(w)\in h^{\tilde X}_wH\o$ and $Y(w)\in V_wHM$ such that $h(w)+Y(w) \in \tilde E^u(w)$. We then have for any $t\in\R,$ $\tilde T^t(h(w))=H_{\tilde X}(\tilde T^t(Y(w)))$. Note $Y$ the vector field defined along $\tilde \ph.w$ by $Y(\tilde \ph^t(w))=\tilde T^t(Y(w))$. Since 
$$\tilde T^t(h(w))=H_{\tilde X}(Y)(\tilde \ph^t(w))=-[\tilde X,Y](\tilde \ph^t(w)),$$
we get in a good affine chart at $w$,
$$\begin{array}{llllr}
\tilde T^t(h(w))&=&-\ (m(w))^{\frac{1}{2}} (L_{\tilde X} m^{\frac{1}{2}})(\tilde \ph^t(w))\ .\ T^t_e(Y(w))\\\\

&& -\ (m(w)m(\tilde \ph^t(w)))^{\frac{1}{2}} \ .\ T_e^t(h(w)).\\ 
& \end{array}$$
We then get the result since the first term is killed by $d\pi$.
\end{proof}

\begin{figure}[h!]
\begin{center}
\includegraphics[angle=0,width=9cm]{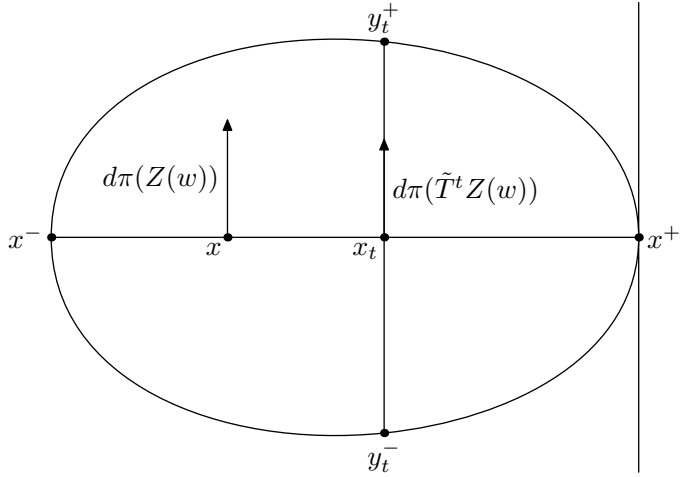} 

\end{center}

\caption{Parallel transport on $\ho$}
\label{oh}
\end{figure}

If $f$ and $g$ are two functions of $t\in\R$, $f(t) \asymp g(t)$ will mean that
$f(t)=O(g(t))$ and $g(t)~=~O(f(t)),$
that is there exists $C>0$ such that
$C^{-1} |f(t)| \leq |g(t)| \leq C |f(t)|.$\\\\
The following proposition gives a link between the parallel transport and the boundary of $\o$.
\begin{prop}\label{transport}
Let  $Z(w)\in \tilde E^{u,s}(w)$. In a good chart at $w$, 
$$\|\tilde T^t(Z(w))\| \asymp \left(\frac{|x_tx^+|^{1/2}}{x_ty_t^+}+\frac{|x_tx^+|^{1/2}}{x_ty_t^-}\right),$$
where $x_t=\pi(\tilde \ph^t(w))$ and $y_t^{\pm}$ are the intersections of the line $x_t+\lambda d\pi(\tilde T^t Z(w))$ with the boundary $\partial \Omega$. (c.f. Figure \ref{oh})
\end{prop}

\begin{proof}
Let choose a good chart at $w$. We have for $Z(w)=Y(w)+h(w) \in \tilde E^{u,s}(w)$,
$$\|\tilde T^t(Z(w))\| \asymp F(d\pi(\tilde T^t(h(w))) \asymp |x_tx^+|^{1/2} F(d\pi(T_e^t(h(w))))$$
from lemma \ref{horver}. But
$$F(d\pi(T^t_e h(w)))=|d\pi(T^t_e h(w))|m^{-1}(d\pi(T^t_e h(w)))=
|d\pi(h(w))| m^{-1}(d\pi(T^t_e h(w))).$$
Since
$$m^{-1}(d\pi(T^t_e h(w)))=\frac{1}{2}\left(\frac{1}{x_ty_t^+}+\frac{1}{x_ty_t^-}\right),$$
we get
$$\|\tilde T^t(Z(w))\| \asymp \left(\frac{|x_tx^+|^{1/2}}{x_ty_t^+}+\frac{|x_tx^+|^{1/2}}{x_ty_t^-}\right).$$
\end{proof}

\subsection{Proof of the Anosov property}\label{proofanosov}

\begin{lemma}\label{equivalents}
In a good chart at $w=(x,[\xi])$ holds the following asymptotic expansion :
$$x_tx^+ = \frac{|xx^+|^2}{m(w)}e^{-2t} + O(e^{-4t}).$$
\end{lemma}
 
\begin{proof}
From the fact that $d_{\o}(x,x_t) = t$, a direct calculation yields
$$\begin{array}{lll}
xx_t=\displaystyle\frac{e^{2t}-1}{\frac{1}{xx^-}+\frac{1}{xx^+}e^{2t}}&=xx^+(1-e^{-2t})(1-\displaystyle\frac{xx^+}{xx^-}e^{-2t} + O(e^{-4t}))\\\\
&=xx^+(1-\left(\displaystyle\frac{x^-x^+}{xx^-}\right)e^{-2t} + O(e^{-4t}))
\end{array}$$
Hence
$$x_tx^+=xx^+-x_tx=x^-x^+\frac{xx^+}{xx^-}e^{-2t} + O(e^{-4t})=\frac{|xx^+|^2}{m(w)}e^{-2t} + O(e^{-4t}).$$
\end{proof}

We can now prove theorem \ref{anosov}.
\begin{proof}[Proof of theorem \ref{anosov}]
Let $E_1^s=\{v\in E^s, \|v\|=1\}$ the set of unit ``stable'' vectors and
$$f : E_1^s \times \R \longrightarrow \R$$
the continuous function defined by
$$f(v,t)= \| T^t(v)\| e^{-t}.$$
Choose $\tilde v \in \tilde E_1^s(x,[\xi]) \subset T_{(x,[\xi])}H\o$ corresponding to $v$. If the chart is adapted to $(x,[\xi])$, just remember that the vector $d\pi(\tilde T^t \tilde{v})$ is orthogonal to $\overrightarrow{x_tx^+}$; hence so are $\overrightarrow{x_ty_t^+}$ and $\overrightarrow{x_ty_t^-}$. Proposition \ref{transport} and lemma \ref{equivalents} yield
$$\| \tilde T^t(\tilde v)\|  \asymp \frac{1}{|x_tx^+|^{1/2}}\left(\frac{x_tx^+}{x_ty_t^+}+\frac{x_tx^+}{x_ty_t^-}\right)
 \asymp e^{t}\left(\frac{x_tx^+}{x_ty_t^+}+\frac{x_tx^+}{x_ty_t^-}\right),$$ 
 From the strict convexity of $\o$, we have $\lim_{t\to\infty} x_tx^+ / x_ty_t^{\pm}=0$ ; we thus deduce
$$\| T^t(v)\|=\| \tilde T^t(\tilde v)\| \ll e^t,$$ so that for any $v\in E_1^s$, 
$$\lim_{t\to\infty} f(v,t)=0.$$
Choose $0<a<1$. Since $E_1^s$ is compact and $f$ continuous, we can find a time $t_a> 0$, such that for any $t\geq t_a$ and any $v\in E_1^s$, 
$f(v,t)\leq a$, meaning that for any $v\in E^s$, 
$$\| T^t(v)\| \leq a\|v\|e^{t}.$$
Let $M_a = \max \{ \|T^t\|, \ 0\leq t\leq t_a \}$. Iterating the precedent inequality, we get for $t$ large enough and any $v\in E^s$,
$$\| T^t(v)\| \leq a\|T^{t-t_a}(v)\|e^{t_a} \leq \cdots \leq a^{[t/t_a]} \|T^{t-[t/t_a]t_a}(v)\|e^{[t/t_a]t_a}  \leq M_a a^{t/t_a} e^{t} = M_a e^{(1-\alpha)t}$$
with $\alpha=-\log(a)/t_a>0.$ Hence, from equation (\ref{dd}), for any $v\in E^s$,
$$\|d\ph^t(v)\| \leq M_a e^{-\alpha t}.$$
Reversing the time and using $J^X$, we get the result for $v\in E^u$, which completes the proof.
\end{proof}


\vskip 1cm
\section{Lyapunov exponents}\label{sectionlyap}
\subsection{Generalities}

\begin{defi}\label{defilyapunov}
Let $\ph=(\ph^t)$ be a $C^1$ flow on a manifold $W$.
The point $w\in W$ (or his orbit $\ph.w$) is regular if there exists a $\ph^t$-invariant decomposition
$$TW=\R.X + \oplus_{i=1}^{p} E_i$$
along $\ph.w$ and real numbers 
$$\ch_1(w) \leq \cdots \leq \ch_{p}(w),$$
called Lyapunov exponents, such that, for any vector $v_i\in  E_i \backslash \{0\}$,
$$\lim_{t\to\pm\infty} \frac{1}{t} \log \|d\ph^t(v_i)\| = \ch_i(w).$$
\end{defi}

\begin{thm}[Osedelec's ergodic multiplicative theorem, \cite{osedelec}]\label{fullmeasure}
For any $\ph^t$-invariant measure, the set $\Lambda$ of regular points is of full measure.
\end{thm}

Let us come back to our case, and pick a regular point $w\in \Lambda \subset HM$. Obviously, the Lyapunov decomposition in definition \ref{defilyapunov} will be a subdecomposition of the Anosov decomposition, that is
$$THM = R.X \oplus E^s \oplus E^u = R.X \oplus (\oplus_{i=1}^{p} E_i^s) \oplus (\oplus_{j=1}^{q}  E_j^u).$$
The positive Lyapunov exponents will come from the unstable distribution and the negative from the stable. The following proposition relates the Lyapunov exponents and the parallel transport. Together with proposition \ref{transport}, we get a link between the Lyapunov exponents and the shape of the boundary $\partial\o$.

\begin{prop}\label{core}
Let $w\in\Lambda$ be a regular point. The Lyapunov distribution is given by
$$THM = R.X \oplus \left(\oplus_{i=1}^{p} (E_i^s \oplus E_i^u)\right),$$
with $E_i^s=J^X(E_i^u)$. Furthermore, the corresponding Lyapunov exponents are given by 
$$\ch_i^{\pm}(w) = \pm 1 + \eta_i(w)$$
where 
$$-1 < \eta_1(w) <\cdots < \eta_p(w) < 1$$
are the Lyapunov exponents of the parallel transport $T^t$ at $w$.
\end{prop}
\begin{proof}
Choose $Z_i^u(w)\in E_i^u(w)$ corresponding to the Lyapunov exponent $\ch_i^+(w)$. Then, from equations $(\ref{dd})$, 
$$\ch^+_i(w)=\c(w,Z_i^u(w))=\lim_{t\to\infty}  \frac{1}{t} \log \|d\ph^t(Z_i^u(w))\|= 1 + \lim_{t\to\infty}  \frac{1}{t} \log \|T^t(Z_i^u(w))\| = 1 + \eta_i(w).$$
More over,
$$\lim_{t\to\infty}  \frac{1}{t} \log \|d\ph^t(J^X(Z_i^u(w)))\|= -1 +\lim_{t\to\infty}  \frac{1}{t} \log \|T^t(J^X(Z_i^u(w)))\| = -1 + \eta_i(w),$$
from the remark following lemma \ref{comparison}.
\end{proof}

\subsection{Shape of the boundary}
Here we precise the relation between the Lyapunov exponents and the boundary $\partial\o$. For this we come back to the function
$$g(t,Z(w))=\left(\frac{|x_tx^+|^{1/2}}{x_ty_t^+}+\frac{|x_tx^+|^{1/2}}{x_ty_t^-}\right),$$
which appears in the proposition \ref{transport}.\\
We know from \cite{benoistcv1} that, in our context of a strictly convex set, the metric space $(\o,d_{\o})$ is Gromov-hyperbolic. Then proposition 1.8 of \cite{benoistqs} tells us that
$$x_ty_t^+ \asymp x_ty_t^-$$
since the points $y_t^+,x^-,y_t^-,x^+$ is a harmonic ``quadruplet'' (see \cite{benoistqs} for (here not relevant) details). Thus, 
\begin{equation}\label{gx}
g(t,Z(w)) \asymp \frac{|x_tx^+|^{1/2}}{x_ty_t^+}.
\end{equation}
Assume $w$ is a regular point, choose $Z_i(w)\in \tilde E^{u,s}_i(w)$ and look at the asymptotic exponential behavior of the function $g(t,Z_i(w))$ : we have
$$\lim_{t\to\infty} \frac{1}{t} \log g(t,Z_i(w)) = \eta_i(w),$$
that is, for any $\epsilon>0$,
\begin{equation}\label{encadrement}
e^{(\eta_i(w)-\epsilon)t}\leq g(t) \leq e^{(\eta_i(w)+\epsilon)t} 
\end{equation}
 for $t$ large enough.\\\\
What does this mean on the boundary ? Choose distinct points $x^+,x^-$ on $\partial\o$ and an affine chart in which $T_{x^+}\partial\o$ and $T_{x^-}\partial\o$ are parallel. We choose the euclidean metric such that the segment $[x^+x^-]$ could be identified with $[0,1]$ : a point $x\in (x^+x^-)$ is thus seen as a real in $(0,1)$, and we have $x^+=0,\ x^-=1.$\\
Given a vector $v \in T_{x^+}\partial\o$, we look at the section of $\o$ by the plane $\vect \{ v,\overrightarrow{x^+x^-}\}$, and call $y^{\pm}(v,x)$ the distance from $x$ to the boundary points $y^{\pm}(x)$, intersections of $\partial\o$ and the line $\{x \pm \l v\}_{\l >0}$ (see figure \ref{figbord}).

\begin{figure}[h!]
\begin{center}
\includegraphics[angle=0,width=9cm]{touteslesfigures.6}
\end{center}
\caption{}
\label{figbord}
\end{figure}

We have the following
\begin{prop}\label{bord}
Assume the line $(x^+x^-)$ is the projection of a regular orbit of the flow, with Lyapunov exponents $\ch_i^{\pm}=\pm 1 + \eta_i,\ i=1\cdots p$. 
Then there exists a decomposition of the tangent space 
$$T_{x^+}\partial\o=\oplus_i^p H_i(x^+),$$
such that, for any $v_i \in H_i(x^+)$,
$$y^{\pm}(v_i,x) \asymp x^{(1+\eta_i)/2}$$
for small $x$.
\end{prop}

\begin{proof}
We first use lemma \ref{equivalents} with $w=(x,[\xi]),\ x$ being the middle point of  the segment $[x^+x^-]$ and $\xi = \overrightarrow{xx^+}$. We have thus
$$xx^+=xx^-=m(w)=\frac{1}{2},$$ which gives
$$x_t = |x_tx^+| = \frac{1}{2} e^{-2t} (1 + o(1)).$$
 Hence
\begin{equation}\label{xt} 
t = \log (x_t^{-1/2}) + O(1).
\end{equation}
Note $F_i = d\pi(\tilde E^s_i)$, $H_i(x^+)=x^+ + F_i$, and pick $v_i\in H_i(x^+)$. Note 
$y_t^{\pm}=y^{\pm}(v_i,x_t)$. From (\ref{gx}) and (\ref{encadrement}), there exists $0<C<1$ such that
$$C^{-1} e^{(\eta_i(w)-\epsilon)t}x_t^{-1/2}\leq \frac{1}{y^{\pm}_t} \leq Ce^{(\eta_i(w)+\epsilon)t}x_t^{-1/2};$$
hence, using (\ref{xt}),
$$D^{-1} x_t^{-(\eta_i(w)+1)/2}\leq y^{\pm}_t \leq D x_t^{-(\eta_i(w)+1)/2},$$
for a constant $0<D<1$.\\
\end{proof}

Note the following : when $\o$ is an ellipsoid,  every point is regular and all the $\eta_i$ are $0$ ; $-1$ and $1$ are the only Lyapunov exponents.
In the next section, we see that if $\o$ is not an ellipsoid, then the Lyapunov exponents vary from a point to another. But we do not know yet if there can be various positive Lyapunov exponents at the same point $w$. (And our paper will let this question unsolved...)

\subsection{Lyapunov exponents of a periodic orbit}\label{periodic}
Every periodic orbit on $HM$ corresponds to a unique non trivial element $\g$ of the group $\Gamma$. As we know from \cite{benoistcv1}, every such element is biproximal, that is : if $(\l_i)_{1\leq i\leq n}$ are its (non-necessary distinct) eigenvalues ordered as $|\l_1| \geq |\l_2| \cdots \geq |\l_{n+1}|$, then $|\l_1| > |\l_2|$ and $|\l_{n+1}| < |\l_n|$. The length of this periodic orbit is given by
$$l_{\g}=\frac{1}{2}(\log |\l_1| - \log |\l_{n+1}|).$$
Let us do the study in dimension 2. Take an element $\gamma\in\Gamma$ conjugated to the matrix
$$\left(\begin{array}{ccc} 
\l_1 & 0 & 0 \\
0 & \l_2 & 0 \\
 0 & 0 & \l_3 \\
\end{array}\right) \in SL_3(\R)
$$ with  $\l_i\in \C,\ |\l_1| > |\l_2| > |\l_3|$.
The line $(\g^-\g^+)$ is his axis and $\g^0$ his third fixed point. We look at the picture in the affine chart given by the plane $\{x_1+x_3=0\}\subset\R^3$, with the following coordinates:
$$ \g^-=[0:0:1],\ \g^+=[1:0:0],\ \g^0=[0:1:0].$$

This is a good affine chart for the periodic orbit we are looking to. Choose a point $x\in (\g^-\g^+)$ with coordinates $[a_0:0:1-a_0]$ where $a_0\in (0,1)$ and let $w=(x,[\overrightarrow{\g^-\g^+}])$. The point $x_n=\g^n.x$ is given by $$x_n=[a_n:0:1-a_n],$$with $$a_{n+1} = \frac{\l_1 a_n}{\l_1 a_n + \l_2 (1-a_n)}.$$
Now, we look at a vector $v=\overrightarrow{xm} \in \overrightarrow{\g^-\g^+}^{\perp}$ with $m=[a_0:b_0:1-a_0],\ b_0\in\R$. Let  $m_n=\g^n.m = [a_n:b_n:1-a_n]$, $v_n=\overrightarrow{x_nm_n}$, so that $|v_n|=|b_n|$. Then $(b_n)$ is given by

$$b_{n+1}=\frac{\l_2 b_n}{\l_1 a_n + \l_2 (1-a_n)}= \frac{\l_2}{\l_1} \frac{a_{n+1}}{a_n} b_n,$$
which leads to 
$$b_n=\left(\frac{\l_2}{\l_1}\right)^n \frac{b_0}{a_0} a_n.$$
Since $\lim_{n\to\infty} a_n=1$, we get
$$b_n \asymp \left(\frac{\l_2}{\l_1}\right)^n.$$

Let $Z(w)\in T_w H\o$ such that $d\pi(Z(w))=v$. Since $\g$ is an isometry for $F$, we have, with the notations of proposition \ref{transport},
$$\begin{array}{llll}
1\asymp F(x,v)=F(x_n,v_n) & \asymp & \left|\displaystyle\frac{\l_2}{\l_1}\right|^n \displaystyle\frac{1}{|x_n\g^+|^{1/2}} \left(\displaystyle \frac{|x_n\g^+|^{1/2}}{x_ny_n^+}+\displaystyle\frac{|x_n\g^+|^{1/2}}{x_ny_n^-}\right)\\\\
& \asymp & \left|\displaystyle\frac{\l_2}{\l_1}\right|^n  e^{n l_{\g}} \|T^{n l_{\g}}(Z(w))\|,\\
\end{array}$$
by using lemma \ref{equivalents}. Thus
$$\|T^{n l_{\g}}(Z(w))\| \asymp \left|\frac{\l_1}{\l_2}\right|^n  e^{-n l_{\g}}$$
and
$$\lim_{t\to +\infty} \frac{1}{t} \log \|T^{t}(Z(w))\| = \lim_{n\to\infty} \frac{1}{nl_{\g}} \log \|T^{n l_{\g}}(Z(w))\| = -1+2\ \frac{\log \left| \l_1 /\l_2 \right|}{\log \left|\l_1/\l_3\right|}.$$

All this can be generalized to any dimension by sectioning the convex set, so that we get the following result, already known by Yves Benoist \cite{benoistcv1} (but stated in another form and context).

\begin{prop}
The Lyapunov exponents $(\eta_i(\g))$ of the parallel transport along a periodic orbit corresponding to $\g\in\Gamma$ are given by
$$\eta_i(\g) = -1 + 2\ \frac{\log \l_0-\log \l_i}{\log \l_{0} -\log \l_{p+1}},\ i=1\cdots p,$$
where $\l_0 > \l_1 > \cdots > \l_{p} > \l_{p+1}$ denote the moduli of the eigenvalues of $\g$. The corresponding Lyapunov exponents are given by
$$\chi^+_i(\g) = 2\ \frac{\log \l_0-\log \l_i}{\log \l_{0} -\log \l_{p+1}},\ i=1 \cdots p,$$
$$\chi^-_i(\g) = -2 + 2\ \frac{\log \l_0-\log \l_i}{\log \l_{0} -\log \l_{p+1}},\ i=1\cdots p.$$
\end{prop}

Note that in the case of a hyperbolic structure, $p=1$ and $\l_1=1$, so that $\eta_1=0$. In fact, we can find a Riemannian metric $\|.\|$ on $HM$ for which the parallel transport is an isometry. In the other cases, the proposition proves that it is not possible anymore.


\vskip 1cm
\section{Symmetric considerations}\label{symmetry}

At the end of this section will be proved the upper bound in the main theorem \ref{majeur}. Let us recall the point that we reached here.
We know that the Lyapunov exponents can be written 
$$\ch_i^{\pm}=\pm 1 + \eta_i,\ i=1\cdots p.$$
Thus $$\ch^+=\sum_{i=1}^p \dim E_i\ \ch_i^+ = (n-1) + \eta,$$
where $\eta = \sum_{i=1}^p \dim E_i\ \eta_i$, so that we get from Ruelle inequality 
$$h_{top} \leq (n-1) + \int_{HM} \eta\ d\mu_{BM}.$$
In this section we prove that
$$ \int_{HM} \eta\ d\mu_{BM} = 0.$$
Since the measure $\mu$ is ergodic and $\eta$ is $\ph^t$-invariant, this is equivalent to the fact that $\eta = 0$ almost everywhere for $\mu_{BM}$. However, as we saw in the section \ref{periodic}, $\eta$ is not identically $0$ on $\Lambda$. Another question is to know if every point is regular, that is $\Lambda=HM$.\\

For that, we introduce the map 
$$\begin{array}{lclc} 
\sigma :& HM & \longrightarrow & HM\\ 
& w=(x,[\xi]) & \longmapsto & (x,[-\xi]). 
\end{array}$$ 
which is a $C^{\infty}$ involutive diffeomorphism and we look at symmetries with respect to $\sigma$.

\subsection{Symmetric sets and functions} 
We say that 
\begin{itemize}
\item a subset $A$ of $HM$ is symmetric if $\sigma A=A$;
\item a function $f:A\rightarrow \R$ defined on a symmetric set $A$ is symmetric if $f\circ\sigma=f$, antisymmetric if $f\circ\sigma=-f$.
\end{itemize}

\begin{prop}
\begin{itemize}
\item The application $\sigma$ exchanges the stable and unstable foliations. \\
\item The set $\Lambda$ of regular points is a symmetric set and $d\sigma$ preserves the Lyapunov decomposition by sending $E^s_i(w)$ to $E^u_i(\sigma(w))$, for any $w\in\Lambda$. \\
\item The function $\eta : \Lambda \longrightarrow \R$ is antisymmetric.\\
\end{itemize}
\end{prop}
\begin{proof}
The differential $d\sigma$ of $\sigma$ is an isomorphism and there exists $C>0$ such that
$$ C^{-1} \leq \| d\sigma \| \leq C.$$
Furthermore, we have
$$\ph^{-t}=\sigma\circ \ph^t\circ \sigma,$$ and then for any $w\in HM$ and $Z^u(w)\in E^u(w)$,
$$ \| d_{\sigma(w)}\ph^t( d_w\sigma(Z^u(w))) \| =  \|  d_{\ph^{-t}(w)}\sigma (d_w\ph^{-t}(Z^u(w))) \| \asymp \|d_w\ph^{-t}(Z^u(w))) \|,$$
which proves that 
$$d_w\sigma (E^u(w))=E^s(\sigma(w)),$$
and also
$$d_w\sigma (E^s(w))=E^u(\sigma(w));$$
that is $\sigma$ exchanges the stable and the unstable foliations.
Moreover, if $w\in \Lambda$, then for $Z(w) \in T_wHM$,
$$\lim_{t\to\infty} \frac{1}{t} \log \| d_w \ph^{-t} (Z(w))\|= - \lim_{t\to\infty} \frac{1}{t} \log \| d_w \ph^{t} (Z(w))\|=-\chi(w,Z(w)).$$
thus
$$-\chi(w,Z(w)) = \lim_{t\to\infty} \frac{1}{t} \log \| d_w \ph^{-t} (Z(w))\| = 
\lim_{t\to\infty} \frac{1}{t} \log \| d_{\sigma(w)}\ph^t( d_w\sigma(Z(w))) \| =  \chi(\sigma(w),d_w\sigma(Z(w))),$$
which proves that $\sigma(w)$ is also regular, hence $\Lambda$ is symmetric.
We also get the decomposition
$$T_{\sigma(w)}HM= R.X(\sigma(w)) \oplus \left(\oplus_{i}^{p} (E^{s}_i(\sigma(w)) \oplus E^{u}_i(\sigma(w)))\right)$$
with
$$E^{s}_i(\sigma(w))= d\sigma(E_i^u(w)),\ E^{u}_i(\sigma(w))= d\sigma(E_i^s(w)).$$
Furthermore,
\begin{equation}\label{lyapu-}
\ch^+_{i}(\sigma(w))=-\ch_{p+1-i}^-(w),
\end{equation}
so that
$$\eta_{i}(\sigma(w))=-\eta_{p+1-i}(w).$$
We finally get
$$\eta(\sigma(w))=\sum_{i=1}^p \dim E_i(\sigma(w))\ \eta_i(\sigma(w))=-\eta(w).$$
\end{proof}

\subsection{Symmetric measures}

Let $\mathcal{M}$ be the set of invariant probability measures for $\ph$. This is a non empty convex set whose extrem points are ergodic measures, that is measures such that every invariant Borel set has either full or zero measure.\\
For any measure $\mu\in\M$, we can consider $\sigma\star\mu$ ; we say that $\mu$ is symmetric if $\sigma\star\mu=\mu$. $\M^s$ will denote the set of symmetric measures.\\
Any function $f$ on $HM$ can be written as the sum of a symmetric and an antisymmetric function, that is $f=g+h$ with $g\sigma=g$, $h\sigma=-h$. Furthermore, if $\nu\in\M^s$, then
$\int h d\nu=0$, that is $\int f\ d\nu = \int g \ d\nu$.
\begin{lemma}
$\M^s$ is a weakly closed convex subset of $\M$.
\end{lemma}
\begin{proof}
It is obviously convex. To prove that is closed, let $(\mu_n)$ be any sequence in $\M^s$ converging to $\mu\in\M$ in the weak sense, that is, for any $f\in C^0(HM)$, 
$$\lim_{n\to\infty} \int f\ d\mu_n=\int f\ d\mu.$$
Since the $\mu_n$ are symmetric, $(\sigma\star\mu_n)$ also converge to $\mu$. But, with the decomposition $f=g+h$, 
$$\int h \ d\mu_n =0,\ \int f \ d\mu_n = \int g \ d\mu_n,$$ 
so we have, for any $f\in C^0(HM)$,
$$\int f\ d(\sigma\star\mu) = \int g\ d\mu - \int h \ d\mu = \lim_{n\to\infty} \int g\ d\mu_n - \int h\ d\mu_n = \int f\ d\mu,$$
which gives $\mu=\sigma\star\mu$.
\end{proof}

\begin{prop}
The Bowen-Margulis measure on $HM$ is symmetric.
\end{prop}
\begin{proof}
We only have to recall the construction of the measure made by Bowen in \cite{bowen1} (c.f. \cite{kh}):  it is obtained by taking a weak limit of the family of symmetric measures given by
$$\mu_{t,\epsilon}=\frac{1}{l(t,\epsilon)} \sum_{\g\in \textrm{Per}(t,\epsilon)} \delta_{\g}.$$
Here, $\textrm{Per}(t,\epsilon)$ denotes the set of closed orbits $\g$ with a period in $(t,t+\epsilon)$ and $\delta_{\g}$ the Lebesgue measure on $\g$ ; $l(t,\epsilon)$ is the sum of the lengths (minimal periods) of the $\g\in\textrm{Per}(t,\epsilon)$. Note that if $\g$ is in $\textrm{Per}(t,\epsilon)$, so is the symmetric orbit $\sigma\star\g$. Then every $\mu_{t,\epsilon}$ is symmetric, and any weak limit also. 
\end{proof}

All this stuff proves the first part of theorem \ref{majeur}, that is
\begin{prop}
Let $\ph$ be the geodesic flow on the Hilbert metric on a compact strictly convex projective manifold $M$ of dimension $n$. Its topological entropy $h_{top}(\ph)$ satisfies the inequality
$$h_{top}(\ph) \leq (n-1).$$ 
\end{prop}
\begin{proof}
Since $\mu_{BM}$ is symmetric and $\eta$ antisymmetric, we have $\int \eta\ d\mu_{BM} =0$, which yields 
$$h_{top}(\ph)=h_{\mu_{BM}} \leq (n-1).$$
\end{proof}


\vskip 1cm

\section{Invariant measures and the equality case}\label{equality}

\subsection{The equality case}
Here we deal with the equality case in theorem \ref{majeur}. This is closely related to the equality case in the Ruelle inequality (\ref{ruelle}), that is : for which measures $\mu\in \M$ do we have 
$$h_{\mu} = \int \ch^+\ d\mu\ ?$$
Ledrappier and Young answered this question in the first part of \cite{ly} :
\begin{thm}[\cite{ly}, Theorem A]\label{lyequality}
Let $\ph^t :W\longrightarrow W$ be a $C^{1+\epsilon}$ flow on a compact manifold $W$. Then an invariant measure $\mu$  has absolutely continuous conditional measures on unstable manifolds if and only if
$$h_{\mu} = \int_W \ch^+\ d\mu.$$
\end{thm}
(In the original paper, this is proved for $C^2$ diffeomorphisms, but it extends to our case. See \cite{pesinbarreira} for a complete presentation.)

From this theorem we can now prove the
\begin{prop}
$h_{top} = n-1$ if and only if the Hilbert metric is Riemannian.
\end{prop}
\begin{proof}
$h_{top} = n-1$ if and only if $h_{\mu_{BM}}=n-1$, that is the Bowen-Margulis measure satisfies the equality in the Ruelle inequality. But from theorem \ref{lyequality}, it is equivalent to the absolute continuity of its conditional measures on unstable manifolds, that is the absolute continuity of the Margulis measures $\mu^u$ on strong unstable manifolds : recall that Margulis constructed the (future) Bowen-Margulis measure as a local product $\mu^s \times \mu^u \times dt$, where the measures $\mu^s$ and $\mu^u$ were measures on strong stable and strong unstable manifolds with adequate properties (see \cite{margulis}, \cite{margulisbook} or \cite{kh} for more details). It follows from the symmetry of this construction that $\mu^u$ is absolutely continuous if and only if $\mu^s$ is so, that is if and only if $\mu_{BM}$ is absolutely continuous. The proposition \ref{absco} concludes the proof.
\end{proof}

We can add some remarks to this proof and connect it with some well known results in the ergodic theory of hyperbolic systems. In our context of a topologically transitive Anosov flow with dense periodic orbits, we indeed know from \cite{bowenruelle} that there exists only one invariant measure $\mu^+$, called the Sinai-Ruelle-Bowen (SRB) measure, which satisfies the equality in (\ref{ruelle}). This measure is ergodic and characterized by either of the following equivalent facts :
\begin{itemize}
\item $\mu^+$ satisfies the equality in (\ref{ruelle}) ;
\item the conditional measures $(\mu^+)^u$ on unstable manifolds is smooth ;
\item the equality \begin{equation}\label{average} \lim_{T\to\infty} \frac{1}{T} \int_{0}^{T} f(\ph^t(x))\ dt = \int f(x) \ d\mu^+(x), 
\end{equation}
holds for $\l$-almost every point $x\in HM$.\\
\end{itemize}
Reversing the time, we get the SRB measure $\mu^-$ for $\ph^{-t}$, which is equal to $\mu^+$ if and only if one of the two measures is smooth. Roughly speaking, those two measures are the smoothest invariant measures of the system.\\
In the case of a hyperbolic geodesic flow, the Bowen-Margulis and the SRB measures (for $\ph^t$ and $\ph^{-t}$) coincide with the Liouville measure. This is not true anymore when  the Hilbert metric is not Riemannian : we then get three measures which are of interest, each one being singular with respect to the others. The two measures $\mu^+$ and $\mu^-$ are related via $\sigma$ by $\mu^+=\sigma\star\mu^-$ ; hence $\sigma$ is a smooth diffeomorphism of $HM$ which sends the measure $\mu^-$ to the measure $\mu^+$ which is singular with respect to $\mu^-$.\\\\
The function $\eta$ is invariant under the flow and thus is  constant almost everywhere with respect to either of the ergodic measures $\mu_{BM}$, $\mu^-$ and $\mu^+$. We have already seen that $\eta$ was zero $\mu_{BM}$-almost everywhere. But with respect to $\mu^-$ or $\mu^+$, $\eta$ is equal to a constant $\eta_{SRB}<0$ almost everywhere, since we can prove that $h_{\mu^+}=h_{\mu^-}$ and use $h_{\mu^+}=n-1 + \eta_{SRB} < n-1.$\\
This function $\eta$ has a geometric and physical meaning with respect to the parallel transport $T^t$, since
$$\eta(w) 
= \lim_{t\to \pm\infty} \frac{1}{|t|} \log \left|\det \left(T^t|_{E^u(w)}\right)\right| 
= \lim_{t\to \pm\infty} \frac{1}{|t|} \log \left|\det \left(T^t|_{E^s(w)}\right)\right|
= \lim_{t\to \pm\infty} \frac{1}{2|t|} \log \left|\det T^t_w \right|.$$
That is it measures the asymptotic exponential behaviour of the euclidean volume of a basic parallelepiped on the tangent space. From a Bowen-Margulis measure point of view, this volume is thus asymptotically constant, whereas from a SRB measure point of view, it is exponentially decreasing.\\\\
This should also mean something on the boundary by using proposition \ref{bord}. Indeed, we know from \cite{kaimanovich} (see also \cite{coornaert}) that every $\Gamma$-invariant measure on the boundary $\partial\o$ gives rise to an invariant measure $\mu\in\M$. In particular, to the Patterson-Sullivan measure corresponds the Bowen-Margulis measure.\\
Now assume we are plane, that is $n=2$. Then there is only one positive Lyapunov exponent and one $\eta_i$ which is $\eta$. But $\eta=0\ \mu_{BM}-a.e.$, thus proposition \ref{bord} could mean that (very roughly speaking,) the boundary is locally an ellipsoid at almost every point from the Patterson-Sullivan measure point of view.

\subsection{A large lower bound for the entropy}
We conclude this section by giving a lower bound for the topological entropy in terms of regularity of the boundary. When $\o$ is not an ellipsoid, then the boundary is known to be  $C^{\alpha}$ for a certain $\alpha>1$ but the supremum $\alpha_{\o}$ of such $\alpha$'s is stricly less than $2$. Equivalently (see \cite{benoistcv1}), the boundary is $\beta$-convex, for a certain $\beta>2$, that is there exists a constant $C>0$, such that, for any $p,\ p'\in\partial\o$,
$$d_{\R^n}(p',T_p\partial\o) \geq C |pp'|^{\beta},$$
where $d_{\R^n}$ denotes the euclidean distance. The corresponding infimum $\beta_{\o}>2$ satisfies
$$\frac{1}{\beta_{\o}}+\frac{1}{\alpha_{\o}}=1.$$

\begin{prop}
Let $M=\o/\Gamma$, where $\o$ is not an ellipsoid, and assume $\partial\o$ is $\beta$-convex for a $\beta \in (2,+\infty)$. Then 
$$h_{top}(\ph) >  \frac{2}{\beta}(n-1).$$
\end{prop}
\begin{proof}
The $\beta$-convexity of the boundary implies there exists $C>0$ such that, for any $t\geq 0$,
$$x_tx^+ \geq C|y_ty_t^+|^{\beta}.$$
Hence
$$\frac{|x_tx^+|^{1/2}}{y_ty_t^+} \geq D|x_tx^+|^{1/2-1/\beta},$$
for a certain constant $D>0$. Thus any positive Lyapunov exponent $\ch^+_i$ satisfies
$$\ch^+_i = 1 + \eta_i \geq 1 + \lim_{t\to\infty} \frac{1}{t} \log |x_tx^+|^{1/2-1/\beta} =\frac{2}{\beta},$$
from proposition \ref{transport} and lemma \ref{equivalents}. Finally, since $\mu^+$ satisfies the Ruelle entropy formula (\ref{ruelle}), we have
$$h_{top}(\ph) > h_{\mu+} \geq \frac{2}{\beta}(n-1).$$
\end{proof}


\vskip 1cm 

\section{Volume entropy}\label{volume}

On the universal covering $\tilde M$ of a compact Riemannian manifold $(M,g)$, we can consider the volume entropy $h_{vol}(g)$ of $(\tilde M,g)$, which measures the asymptotic exponential growth of  volume of balls in $\tilde M$ :
 $$h_{vol}(g)=\lim_{r\to\infty} \frac{1}{r} \log vol(B(x,r)),$$
where $vol$ denotes the Riemannian volume corresponding to $g$.
In \cite{manning79}, Anthony Manning proved the following result :
\begin{thm}\label{manningvol}
Let $h_{top}$ be the topological entropy of the geodesic flow of $g$ on $HM$. We always have
$$h_{top} \geq h_{vol}(g).$$
Furthermore, if the sectional curvature of $M$ is $<0$ then
$$h_{top} = h_{vol}(g).$$
\end{thm}

In his PhD thesis, Daniel Egloff \cite{egloff} extends this result for some regular Finsler manifolds. Let us see that Manning's proof still works in the special case we are dealing with here, where $F$ is a more ``irregular'' Finsler metric :
\begin{prop}
Let $\ph^t : HM \longrightarrow HM$ be the geodesic flow of the Hilbert metric on the strictly convex projective manifold $M=\o/\Gamma$ and $h_{top}$ denote his topological entropy. Then
$$h_{top} = h_{vol}(d_{\o}).$$ 
\end{prop}
The proof is similar to the one by Manning and we do not reproduce it here. The only point we have to check is the following technical lemma that Manning proved using negative curvature. Here we can compute it directly.
\begin{lemma}
The distance between corresponding points of two geodesics $\sigma, \tau : [0,r]\rightarrow \Omega$ is at most $d_{\Omega}(\sigma(0),\tau(0))+d_{\Omega}(\sigma(0),\tau(0))$.
\end{lemma}
\begin{proof} There are two cases : either $\sigma$ and $\tau$ meet each other or not. Anyway, by joining the point $\sigma(0)$ and $\tau(r)$ with a third geodesic, we see we only have to prove that the distance between two different lines going away from the same point (but not necessary with the same speed) increases.\\
So suppose $c,c' : \R \rightarrow \Omega$ are two lines beginning at the same point $m=c(0)=c'(0)$. Take two pairs of corresponding points $(a,a')=(c(t_1),c'(t_1)),  (b,b')=(c(t_2),c'(t_2))$ with $t_2>t_1\geq 0$. We want to prove that $d_{\Omega}(a,a')<d_{\Omega}(b,b')$. As it is obvious if $t_1=0$, assume $t_1>0$ and note $x,x'$ and $y,y'$ the points on the boundary $\partial\Omega$ of $\Omega$ such that $x,a,a',x'$ and $y,b,b',y'$ are on the same line, in this order. Note also $Y = (mx)\cap (bb')$ and $Y = (mx') \cap (bb')$, so that by convexity of $\Omega$, the six points $Y,y,b,b',y',Y'$ are different and on the same line, in this order. The two lines $(aa')$ and $(bb')$ meet at a certain point that we can send at infinity by an homography. So we can assume the two lines are parallel (c.f. figure \ref{figmanning}).\\
From now we only need Thales' help to see that $$1 > [x,a,a',x] = [Y,b,b',Y'] > [y,b,b',y'],$$ so that $$d_{\Omega}(a,a') = |\log([x,a,a',x])| < |\log([y,b,b',y'])| = d_{\Omega}(b,b').$$ 
\end{proof}

\begin{figure}
\begin{center}
\includegraphics[angle=0,width=9cm]{touteslesfigures.7}
\end{center}
\caption{}
\label{figmanning}
\end{figure}

As a corollary of this proposition and theorem \ref{majeur}, we get corollary \ref{vol}.

\bibliographystyle{plain}
\bibliography{entropies}

\end{document}